\documentclass[12pt,reqno]{amsart}

\usepackage{amssymb}

\theoremstyle{plain}
\newtheorem{theorem}{Theorem}
\newtheorem{lemma}{Lemma}
\newtheorem{corollary}{Corollary}
\newtheorem*{conjecture}{Conjecture}
\newtheorem{problem}{Problem}

\newcommand{\lemlabel}[1]{\label{lem:#1}}
\newcommand{\thmlabel}[1]{\label{thm:#1}}
\newcommand{\corolabel}[1]{\label{coro:#1}}
\newcommand{\problabel}[1]{\label{prob:#1}}

\newcommand{\lemref}[1]{\ref{lem:#1}}
\newcommand{\thmref}[1]{\ref{thm:#1}}
\newcommand{\cororef}[1]{\ref{coro:#1}}
\newcommand{\probref}[1]{\ref{prob:#1}}

\newcommand\UU{\mathcal{U}}
\newcommand\QQ{\mathcal{Q}}
\newcommand\JJ{\mathcal{J}}

\newcommand\iv{^{-1}} 

\title[Strongly Right Alternative Rings and Bol Loops]
{Strongly Right Alternative Rings and Bol Loops}

\author[M.~K.~Kinyon]{Michael~K.~Kinyon}
\address{Department of Mathematical Sciences \\
Indiana University South Bend \\
South Bend, IN 46634 USA}
\email{mkinyon@iusb.edu}
\urladdr{http://www.iusb.edu/\symbol{126}mkinyon}
\author[J.~D.~Phillips]{J.~D.~Phillips}
\address{Department of Mathematics \& Computer Science \\
Wabash College \\
Crawfordsville, IN 47933 USA}
\email{phillipj@wabash.edu}
\urladdr{http://www.wabash.edu/depart/math/faculty.html{\#}Phillips}


\begin{document}

\begin{abstract}
We partially answer two questions of Goodaire by showing that in a finite, strongly right alternative ring, the set of units (if the ring is with unity) is a Bol loop under ring multiplication, and the set of quasiregular elements is a Bol loop under ``circle'' multiplication.
\end{abstract}

\maketitle


A \emph{magma} $(L,\cdot)$ is a set $L$ with a binary operation
$\cdot :L\times L\to L$. 
A magma $(L,\cdot)$ is a \emph{loop} if there is a neutral element $1\in L$
satisfying $1x = x1 = x$, and if all right translations $R(a): x\mapsto xa$
and all left translations $L(a): x\mapsto ax$ are bijections.
If a magma $L$ satisfies the \emph{right Bol identity}
$(xy\cdot z)y = x(yz\cdot y)$, then we will say that $L$ is a \emph{Bol magma}.
The Bol identity can be equivalently expressed in terms of right translations:
$R(y)R(z)R(y) = R(yz\cdot y)$. A \emph{Moufang loop} is a Bol loop which also
satisfies the \emph{flexible law} $x\cdot yx = xy\cdot x$. For background
in loops, the reader is referred to the standard references \cite{Br} and \cite{Pf},
and further details on Bol loops can be found in \cite{R}.

A (nonassociative) ring $(R,+,\cdot)$ is \emph{right alternative} if it
satisfies the \emph{right alternative law} $xy\cdot y = xy^2$. A right alternative
ring $R$ is \emph{strongly right alternative} if $(R,\cdot)$ is a Bol magma.
A right alternative ring of characteristic not $2$ is strongly right alternative
(\cite{ZSSS}, {\S}16.1), but this is not true in general  \cite{G1} \cite{K}. 
An \emph{alternative} ring is a right alternative ring which also satisfies
the \emph{left alternative law} $x\cdot xy = x^2 y$. An alternative ring
is strongly right alternative (\cite{ZSSS}, {\S}2.3), but there are strongly right alternative rings which are not alternative.

In a ring $R$ with unity, a \emph{unit} is an element with a
(not necessarily unique)
two-sided inverse. Let $\UU(R)$ denote the set of all units in
$R$. It is well-known that if $R$ is an associative ring with unity,
then $\UU(R)$ is a group. More generally, if $R$ is an alternative
ring with unity, then $\UU(R)$ is a Moufang loop (\cite{GJM}, {\S}11.5.3).
In a strongly right alternative ring with unity, it is known that if
$\UU(R)$ is closed under multiplication, then $\UU(R)$ is a Bol loop
(\cite{G2}, Lem. 3.1, p. 356). Goodaire posed the following
problem (\cite{G2}, p. 356):

\begin{problem}
\problabel{units}
In a strongly right alternative ring $R$ with unity, is $\UU(R)$ a Bol loop?
\end{problem}

In a ring $R$, define a binary operation $\circ : R\times R\to R$ by
$x \circ y = x + y + x\cdot y$
for $x,y\in R$. An element $x$ is said to be \emph{quasiregular}
with \emph{quasi-inverse} $x'$ if $x\circ x' = x'\circ x = 0$.
Let $\QQ(R)$ denote the set of all quasiregular elements of $R$.
If $R$ is associative, then $(R,\circ)$ is associative, and so
$(\QQ(R),\circ)$ is a group. More
generally, if $R$ is alternative, then $(\QQ(R),\circ)$ is a 
Moufang loop \cite{G0}. For a right alternative ring $R$,
the magma $(R,\circ)$ is right alternative, that is,
$(x\circ y)\circ y = x\circ (y\circ y)$, and if $R$ is strongly
right alternative, then $(R,\circ)$ is a Bol magma with neutral
element $0$ (\cite{G1}, Lem. 3.1). Thus we have the following problem
\cite{G1}.

\begin{problem}
\problabel{circle}
In a strongly right alternative ring $R$, is $(\QQ(R),\circ)$ a Bol loop?
\end{problem}

In this paper, we offer partial affirmative answers to Problems \probref{units} and
\probref{circle}. Our main result, however, is not actually about right
alternative rings. In what follows, let $(L,\cdot)$ denote a Bol magma
with neutral element $1\in L$,
and let $\JJ(L)$ denote the subset of all elements of $L$ with two-sided
inverses, that is, for each $a\in \JJ(L)$, there exists $a\iv\in L$ such
that $aa\iv = a\iv a = 1$. (The uniqueness of $a\iv$ will be addressed in
Lemma \lemref{bijections} below.)

\begin{theorem}
\thmlabel{closed}
If $L$ is a finite Bol magma with a neutral element, then $(\JJ(L),\cdot)$ is a Bol loop.
\end{theorem}

Implicit in the statement of Theorem \thmref{closed} is the assertion that $\JJ(L)$
is closed under multiplication. In fact by the following lemma, the
closure is all that needs to be verified.

\begin{lemma}
\lemlabel{main}
If $L$ is a Bol magma with a neutral element, and if $\JJ(L)$ is closed under multiplication,
then $(\JJ(L),\cdot)$ is a loop.
\end{lemma}

A proof of Lemma \lemref{main} can be extracted from the proof of
(\cite{G2}, Lem. 3.1, p. 356). We will give an independent proof for completeness. 

As an immediate consequence of Theorem \thmref{closed}, we obtain our partial
affirmative answers to Problems \probref{units} and \probref{circle}.

\begin{corollary}
\corolabel{loops}
Let $(R,+,\cdot)$ be a finite, strongly right alternative ring.
\begin{enumerate}
\item If $R$ is with unity, then $(\UU(R),\cdot)$ is a Bol loop.
\item $(\QQ(R),\circ)$ is a Bol loop.
\end{enumerate}
\end{corollary}

We will also use Lemma \lemref{main} to obtain an auxiliary result.

\begin{theorem}
\thmlabel{flex}
If $L$ is a flexible Bol magma with a neutral element,
then $(\JJ(L),\cdot)$ is a Moufang loop.
\end{theorem}

As a consequence, we obtain a new proof of the following; see {\S}11.5.3 of
\cite{GJM} and \cite{G0}.

\begin{corollary}
\corolabel{Moufang}
Let $(R,+,\cdot)$ be an alternative ring.
\begin{enumerate}
\item If $R$ is with unity, then $(\UU(R),\cdot)$ is a Moufang loop.
\item $(\QQ(R),\circ)$ is a Moufang loop.
\end{enumerate}
\end{corollary}

Our investigations were aided by the automated deduction
tool OTTER developed by McCune \cite{MC}.

Recall the definition of powers in a magma: $a^0 := 1$, $a^n := a^{n-1}a$
for all positive integers $n$.

\begin{lemma}
\lemlabel{rpa}
For all nonnegative integers $n$, and all $a\in L$, $R(a)^n = R(a^n)$.
\end{lemma}

\begin{proof}
The proof of this result for Bol loops carries through to the
present setting without change; see \cite{R} or (\cite{Pf}, IV.6.5).
\end{proof}

\begin{lemma}
\lemlabel{LR}
If $a,b,c\in L$ satisfy $ab = ca = 1$, then $b = c$.
\end{lemma}

\begin{proof}
For all $x$, the right Bol identity gives $c(ax\cdot a) = xa$.
Taking $x = b^2$ and using Lemma \lemref{rpa}, we get
$c\cdot ba = b^2 a$. Taking $x = b$, we obtain $ba = 1$.
Thus $c = b^2 a = (ba\cdot b^2) a = b(ab^2 \cdot a)
= b\cdot ba = b$, using Lemma \lemref{rpa} and the right
Bol identity again. The rest is clear.
\end{proof}

\begin{lemma}
\lemlabel{bijections}
For each $a\in \JJ(L)$,
\begin{enumerate}
\item[1.] There is a unique two-sided inverse $a\iv\in \JJ(L)$, 
\item[2.] $R(a)$ is a bijection of $L$ and $R(a)\iv = R(a\iv)$,
\item[3.] $L(a)$ is a bijection of $L$ and $L(a)\iv = R(a)L(a\iv)R(a\iv)$.
\end{enumerate}
\end{lemma}

\begin{proof}
(1) is immediate from Lemma \lemref{LR}. For (2), we compute
\begin{equation}
\label{eq1}
R(a)R(a\iv)R(a) = R(aa\iv \cdot a) = R(a)
\end{equation}
by the right Bol identity. Next, set $c := (a\iv)^2$
and compute
\begin{equation}
\label{eq2} 
R(a)R(c)R(a) = R(ac\cdot a) = I
\end{equation}
since $ac = aa\iv\cdot a\iv$ by Lemma \lemref{rpa}.
Thus
\[
R(a\iv)R(a) = R(a)R(c)R(a)R(a\iv)R(a) = R(a)R(c)R(a) = I
\]
using \eqref{eq2}, \eqref{eq1}, and \eqref{eq2} again. A
similar argument gives $R(a)R(a\iv) = I$.

For (3): The right Bol identity can be written as $L(xy)R(y) = L(y)R(y)L(x)$.
Taking $x = a\iv$, $y = a$, and using $R(a)\iv = R(a\iv)$, we have
\[
L(a)[R(a)L(a\iv)R(a\iv)] = I.
\]
On the other hand, taking $x = a$ and
$y = a\iv$, then $[R(a)L(a\iv)R(a\iv)]L(a) = I$. This completes the proof.
\end{proof}

The proof of part (2) of Lemma \lemref{bijections} is an adaptation to
this setting of the proof of Theorem 3.10(3) of \cite{Ki} (see also Remark 5, p.~51).

We now prove our main lemma.

\begin{proof}[Proof of Lemma \lemref{main}]
If $\JJ(L)$ is closed under multiplication, then Lemma \lemref{bijections}
implies that for each
$a\in \JJ(L)$, the restrictions of $R(a)$ and $L(a)$ to $\JJ(L)$ are
bijections. Therefore $\JJ(L)$ is a loop.
\end{proof}

\begin{lemma}
\lemlabel{left}
For each $a,b\in \JJ(L)$, $L(ab)$ is bijective, and
\[
L(ab)\iv = R(b) R(a) L(b\iv a\iv) R(a\iv) R(b\iv).
\]
In particular, $ab$ has a unique right inverse, that is,
there exists a unique $c\in L$ such that $ab\cdot c = 1$.
\end{lemma}

\begin{proof}
By the Bol identity, $L(ab)R(b) = L(b)R(b)L(a)$,
and so by Lemma \lemref{bijections} and the Bol identity again,
\begin{align*}
L(ab)\iv &= R(b) L(a)\iv R(b\iv) L(b)\iv \\
  &= R(b) R(a) L(a\iv) R(a\iv) R(b\iv) R(b) L(b\iv) R(b\iv) \\
  &= R(b) R(a) L(a\iv) R(a\iv) L(b\iv) R(b\iv) \\
  &= R(b) R(a) L(b\iv a\iv) R(a\iv) R(b\iv)
\end{align*}
The right inverse of $ab$ is 
$c := 1 L(ab)\iv = (b\iv a\iv \cdot ba) a\iv \cdot b\iv$.
\end{proof}

\begin{lemma}
\lemlabel{torsion}
Suppose $a,b\in L$ satisfy $ab = 1$, and suppose there exist distinct
positive integers $m,n$ such that $b^m = b^n$. Then $ba = 1$.
\end{lemma}

\begin{proof}
By Lemma \lemref{rpa}, $R(b)^m = R(b)^n$. Applying this to $a$, we
obtain $R(b)^{m-1} = R(b)^{n-1}$. Hence there exists a positive
integer $k$ such that $R(b)^k = I$, \textit{i.e.}, $b^k = 1$.
Thus $b$ has a right inverse which, by Lemma \lemref{LR}, must be $a$.
\end{proof}

Finally, we turn to the proof of our main result.

\begin{proof}[Proof of Theorem \thmref{closed}]
Fix $a,b\in \JJ(L)$. By Lemma \lemref{left},
$ab$ has a right inverse $c$. Since $L$ is finite, there exist distinct positive
integers $m, n$ such that $c^m = c^n$. By Lemma \lemref{torsion}, $ab\in \JJ(L)$.
By Lemma \lemref{main}, $\JJ(L)$ is a loop.
\end{proof}

As the proofs of Lemma \lemref{torsion} and Theorem \thmref{closed} indicate, the finiteness
of $L$ can obviously be relaxed to a finite torsion assumption: for each $a\in L$,
there exist distinct positive integers $m, n$ (depending on $a$) such
that $a^m = a^n$. 

Next we turn to our auxiliary result.

\begin{proof}[Proof of Theorem \thmref{flex}]
Fix $a,b\in \JJ(L)$, and let $c\in \JJ(L)$ be the right inverse of $ab$
given by Lemma \lemref{left}. Then $ab\cdot (c\cdot ab) = (ab\cdot c)\cdot ab
= ab = ab\cdot 1$. By Lemma \lemref{left}, we may cancel $ab$ on the left
to obtain $c\cdot ab = 1$ as claimed. Thus $ab\in \JJ(L)$, and so by 
Lemma \lemref{main}, $\JJ(L)$ is a loop.
\end{proof}

\begin{proof}[Proof of Corollary \cororef{Moufang}]
Let $(R,+,\cdot)$ be an alternative ring. If $R$ is with unity, then $R$ is
strongly right alternative and $(R,\cdot)$ is flexible (\cite{ZSSS}, {\S}2.3).
Thus $(\UU(R),\cdot)$ is a flexible Bol magma with a neutral element, and so
by Theorem \thmref{flex}, $(\UU(R),\cdot)$ is a Moufang loop. Also, $(R,\circ)$
is a flexible Bol magma with neutral element $0$ (\cite{G0}, \cite{G1}, Lem. 3.1),
and so Theorem \thmref{flex} implies that $(\QQ(R),\circ)$ is a Moufang loop. 
\end{proof}

By applying Lemma \lemref{main}, one may show that $\JJ(L)$ is a loop
for various other special classes of Bol magmas, but lacking examples or
motivation, we omit these. In general, we suspect the following to be the
case.

\begin{conjecture}
There exists a Bol magma $L$ with a neutral element such that $\JJ(L)$
is not closed under multiplication.
\end{conjecture}

Returning to strongly right alternative rings, it is possible that the ring structure allows further weakening of the torsion assumptions on $(\UU(R),\cdot)$ or $(\QQ(R),\circ)$ in Corollary \cororef{loops}.
In any case, we have not been able to find a (necessarily infinite) strongly right alternative ring $R$ such that either $(\UU(R),\cdot)$ or $(\QQ(R),\circ)$ is not closed under multiplication. Thus Problems 1 and 2 in their full generality remain open.


\end{document}